\documentclass{amsart}
\usepackage{amsfonts, amsmath, amssymb, amsthm}
\usepackage{enumerate}
\usepackage{algorithm, algpseudocode}
\usepackage{color}
\usepackage{verbatim}
\usepackage{multirow}
\usepackage{url}
\usepackage{graphicx}
\MakeRobust{\Call}

\theoremstyle{definition}

\newcommand{\N}{\mathbb{N}}
\newcommand{\Z}{\mathbb{Z}}
\newcommand{\Q}{\mathbb{Q}}
\newcommand{\R}{\mathbb{R}}

\title[Solving the CDP via Machine Learning]{Solving the Conjugacy Decision Problem via Machine Learning}
\begin{document}

\author[J. Gryak]{Jonathan Gryak}
\address{Jonathan Gryak, CUNY Graduate Center, PhD Program in Computer Science, City University of New York}
\email{jgryak@gradcenter.cuny.edu}

\author[R. Haralick]{Robert M. Haralick}
\address{Robert M. Haralick, CUNY Graduate Center, PhD Program in Computer Science, City University of New York}
\email{rharalick@gc.cuny.edu}

\author[D. Kahrobaei]{Delaram Kahrobaei}
\address{Delaram Kahrobaei, CUNY Graduate Center, PhD Program in Computer Science and NYCCT, Mathematics Department, City University of New York}
\email{dkahrobaei@gc.cuny.edu}

\begin{abstract}
Machine learning and pattern recognition techniques have been successfully applied to algorithmic problems in free groups. In this paper, we seek to extend these techniques to finitely presented non-free groups, with a particular emphasis on polycyclic and metabelian groups that are of interest to non-commutative cryptography.\\
\paragraph{\nopunct}
As a prototypical example, we utilize supervised learning methods to construct classifiers that can solve the conjugacy decision problem, i.e., determine whether or not a pair of elements from a specified group are conjugate. The accuracies of classifiers created using decision trees, random forests, and $N$-tuple neural network models are evaluated for several non-free groups. The very high accuracy of these classifiers suggests an underlying mathematical relationship with respect to conjugacy in the tested groups.
\end{abstract}

\maketitle
\section{Introduction}
\paragraph{\nopunct}
Group theory has been a rich source of decision problems dating back to Max Dehn, who in 1911 \cite{dehn1911} articulated the word, conjugacy, and isomorphism problems for finitely generated groups. For a given group $G$, the word problem is to determine, for any word $w\in G$, if $w$ is equivalent to the identity of $G$, while the conjugacy problem is to determine, for any two elements $u,v\in G$, if $u$ is conjugate to $v$.\\
\paragraph{\nopunct}
The exploration of solutions to these problems gave rise to combinatorial group theory - the study of groups via their presentations, i.e., sets of generators and relations. Following the advent of computability theory and building upon the undecidability of the word problem in semigroups by Post \cite{post1947recursive}, Novikov proved  in 1955 that the word problem in groups is in general undecidable \cite{novikov1955algorithmic}. Much subsequent work in combinatorial group theory has determined the decidability of the word, conjugacy, and (to a lesser extent) the isomorphism problem for many classes of groups.\\
\paragraph{\nopunct}
Algorithmic problems in finitely presented groups can be studied through various branches of algebra  - representation theory, combinatorial group theory, and geometric group theory. In \cite{haralick2004pattern}, Haralick et al.\ suggested a machine learning approach to solving algorithmic problems in free groups. In this paper we seek to extend these results to non-free groups. As a prototypical example, we will use machine learning techniques to solve the conjugacy decision problem in a variety of groups. Beyond their utilitarian worth, the developed methods provide the computational group theorist a new digital ``sketchpad" with which one can explore the structure of groups and other algebraic objects, and perhaps yielding heretofore unknown mathematical relationships.\\
\paragraph{\nopunct}
We begin in section 2 with discussion of the work of Haralick et al. in free groups and the general framework that they advance to apply machine learning to group-theoretic problems. In section 3 we leverage the additional structure available in non-free groups to the task of feature extraction - the method by which salient information is gleaned from the training data.\\
\paragraph{\nopunct}
Beyond feature extraction, the machine learning system components of model selection, data generation, and performance evaluation are developed in section 4. The learning models of decision trees, random forests, and $N$-tuple neural networks are introduced, along with their attendant parameters. With respect to data generation we recommend three independent data sets used throughout the supervised learning process. While there are many methods of evaluating classifier performance, accuracy will be the primary metric used.\\
\paragraph{\nopunct}
In section 5 we apply our machine learning system to the specific task of solving the conjugacy decision problem.  We chose six different groups for evaluation, including polycyclic, metabelian, and non-solvable groups. We test the performance of all three learning models with various parameters over all test groups, and, in addition to overall classifier accuracy, provide per class and word length analysis of the performance of the best classifier for each group. We conclude in section 6 with a discussion of our results, as well as outlining the subsequent task of exposing the underlying mathematical relationship intimated by our positive results.\\
\paragraph{\nopunct}
This paper is based in part on the Ph.D. dissertation of Jonathan Gryak at the City University of New York, 2017, written under the supervision of Gryak's advisor Delaram Kahrobaei.
\section{Related Work}
In \cite{haralick2004pattern}, Haralick et al.\ posited that pattern recognition techniques are an appropriate methodology for solving problems in combinatorial group theory. To demonstrate, they constructed a machine learning system for discovering effective heuristics for the Whitehead automorphism problem, a search problem in free groups that uses the successive application of the namesake automorphisms to reduce a word to its minimal length.\\
\paragraph{\nopunct}
As mentioned in \cite{haralick2004pattern}, every machine learning system must contend with the following tasks: data generation, feature extraction, model selection, and evaluation. Once the system is constructed, analysis of the system's performance can yield insights into the nature of the problem at hand, and potentially be used to improve upon it. In the following sections we will delve into each of these aforementioned tasks, showing in the process how these techniques can be extended from free groups to finitely presented groups, and ultimately be adapted to solving the conjugacy decision problem The primary difference in the construction of machine learning systems for free and not-free groups is in feature extraction, which is the focus of the next section.
\section{Feature Extraction in Non-Free Groups}
\paragraph{\nopunct}
One of the most important aspects of creating a machine learning system is the process of feature extraction, the means by which relevant information is distilled from the raw dataset and presented to the learning algorithm.  If the raw dataset is unstructured, subsets of data may first be aggregated into units of observation, from which the features will be extracted. Some datasets may come with an intrinsic structure, such as that of a list, a matrix, or a string of text. Regardless of the data's inherent structure, the ability to extract features from the underlying data that provide information relevant to the learning process requires domain-specific knowledge.\\
\paragraph{\nopunct}
Finitely presented groups, in addition to their representation as generators and relators, have a combinatorial structure that is manifested by their Cayley graphs. A \emph{Cayley graph} is a rooted, labeled digraph, with a vertex for every element in the group and each edge labeled by a generator or an inverse generator. The root of the graph is the identity element. If the group is infinite then so is its Cayley graph. The graph is connected, and the label of every path from the root to a vertex is a word representing a group element. Circuits from the root represent words that are equivalent to the identity and are therefore in the normal closure of the set of relators.\\
\paragraph{\nopunct}
The Cayley graph also enables groups to be considered as metric spaces. Let $G$ be a finitely generated group with generating set $X$ and $\phi:F(X)\rightarrow G$ the canonical epimorphism. Given a word $w$ over the alphabet $X$, let $|w|$ be the \emph{length} of $w$. For $g\in G$, the \emph{geodesic length} of $g$ over $X$ is then defined as
$$
l_X(g)=\min\{|w|\mid w\in F(X), \phi(w)=g\}.
$$
The geodesic length of $g$ corresponds to the shortest path in the Cayley graph whose label represents $g$. If every edge in the Cayley graph is given unit length, then $l_X(g)$ corresponds to the number of edges in the shortest path labeled by $w$. 
Given words $u$ and $v$ representing elements $g$ and $h$ respectively, we can now define the \emph{word metric} $d_X(g,h)=l_X(g^{-1}h)$ that satisfies the axioms required of a metric function. Note that as the notation implies, $d_X(g,h)$ is dependent upon the choice of the generating set $X$. Word length and the word metric provide useful means by which we can associate numerical values with group elements.\\
\paragraph{\nopunct}
The Cayley graph is but one way that graphs can be used to extract numerical information concerning group elements and their word representations. In \cite{haralick2004pattern}, Haralick et al.\ introduced a directed variant of the Whitehead graph that allows one to assign numerical values to subwords and subsequently form feature vectors using these values. Let $F(X)$ be a free group over the alphabet $X$, and $w\in F(X)$. The \emph{labeled Whitehead graph} $\Gamma_W(w)=(V,E)$ for $w$ is an undirected, weighted graph, with the set of vertices $V$ equal to the set $X\cup X^{-1}$, and for every $x_i,x_j\in X\cup X^{-1}$, an edge $(x_i,x_j)$ is added to $E$ if $x_ix_j^{-1}$ or  $x_jx_i^{-1}$ occurs in the cyclic form of $w$. Every edge  $(x_i,x_j)\in E$ is assigned a weight $C(w,x_ix_j)$, which is equal to the number of times the subwords $x_ix_j^{-1}$ or  $x_jx_i^{-1}$ occur in $w$.\\
\paragraph{\nopunct}
These graphs can be generalized to count the occurrence of any pattern of letters within a word. Let $F(X)$ be a free group over $X$, with $\epsilon$ denoting the empty word. For a fixed word $w\in F(X)$ and a finite set of words $U=\{u_1,\ldots,u_k\mid u_j\in F(X)\}$, let us define a weighted directed graph $\Gamma(w)=(V,E)$ for $w$ as follows. The set of vertices $V$ is equal to $X\cup X^{-1}$ as before. For any $x,y\in X\cup X^{-1}$ and $u_j\in U$, we form a directed edge $(x,y)\in E$, labeled by $xu_jy$ and assigned the weight $C(w,xu_jy)$, which is equal to the number of times the reduced subword $xu_jy$ occurs in $w$. Note that unlike in the previous case of the Whitehead graph, we are not considering $w$ cyclically. If $xu_jy$ is equivalent to the empty word in $F(X)$ then no edge is drawn, and $C(w,\epsilon)$ is defined to be 0. If $xu_jy$ is a single letter $a\in X\cup X^{-1}$, then an unlabeled loop $(a,a)$ is added to $E$ and assigned the weight $C(w,a)$. In this framework, every pair consisting of a set of words $U$ and a set of \emph{counting functions} $C(w,xu_jy):F(X)\rightarrow \N$ corresponds to a subgraph of $\Gamma(w)$, which we call a \emph{counting subgraph}.\\
\paragraph{\nopunct}
We can now use these counting subgraphs and their attendant counting functions to extract features from finitely presented groups. Let $G$ be given by the presentation $\langle X|R\rangle$, with $\phi$ the canonical epimorphism, and $U=\{u_1,\ldots,u_k\}$ be a set of words in $F(X)$.  Let $xu_jy$ and $w$ be the geodesic words representing the group elements $\phi(xu_jy)$ (for $1\leq j\leq k$) and $\phi(w)$, and  $C(w,xu_1y),\ldots,C(w,xu_ky)$ be counting functions as defined above. When normalized by the word length $|w|$, these give rise to a real-valued feature vector $v\in\R^k$:
$$
v=\frac{1}{|w|}\langle C(w,xu_1y),\ldots,C(w,xu_ky)\rangle.
$$
\paragraph{\nopunct}
Another means of extracting features from finitely presented groups is via their normal forms. A \emph{normal form} for elements of a group can in general be construed as a unique and most concise representation of each element in the group. For free groups the standard normal form of an element is its reduced word representation. Normal forms need not be words; they can be numbers, sequences, or some other formal representation. Note that in some contexts the uniqueness of normal forms may be relaxed.\\ 
\paragraph{\nopunct}
Every finitely presented group has at least one normal form, as one can impose a total ordering such as shortlex \cite{HEO05} and use the least word under that ordering to represent each group element. However, existence does not entail that the normal form can be calculated efficiently. For some finitely presented groups the Knuth-Bendix algorithm \cite{knuthbendix} can be used to create a confluent, terminating rewriting system with respect to the generating set $X$, resulting in an efficiently calculable normal form for the group. For the class of polycyclic groups the collection algorithm can be used to reduce words to their normal form, and is generally efficient in practice \cite{gryak2016status}.\\

\paragraph{\nopunct}
We now have at our disposal a bevy of combinatorial and geometric machinery to extract information concerning group elements. In the definitions below, let $G$ be given by the finite presentation $\langle X\mid R\rangle$ with $|X|=N$, let $g$ be an element of $G$, and let $w=w_1\ldots w_m$ be a word representation of $g$ over $X$. Let $Y=X\cup X^{-1}$. We first consider feature vectors that are applicable to finitely presented groups that possess an efficiently calculable normal form:
\label{mlfv}
\begin{itemize}
\item $n_0$ (\emph{Normal Form}) - Let $G$ be a finitely presented group possessing an efficiently calculable normal form. If $w$ is a word in normal form, then $w$ is of the form
$$
y_1^{e_1}\cdots y_N^{e_N}
$$
with $y_i\in Y$ and $e_i\in\Z$. The feature vector $n_0$ is then
$$
n_0=\langle e_1,\ldots,e_N\rangle.
$$
\item $n_1$ (\emph{Weighted Normal Form}) - The feature vector $n_1$ is the same as $n_0$ above, except it is weighted by the word length $|w|$:
$$
n_1=\frac{1}{|w|}\langle e_1,\ldots,e_N\rangle.
$$
\end{itemize}
We now consider features for finitely presented groups that do not require a normal form:
\begin{itemize}
\item $f_0$ (\emph{Generator Count}) - Let the generator set $X$ be given a fixed order and let  $x_i\in X$ be the $i$th generator. The counting function $C(w,x_i)=|\{w_j\mid w_j=x_i\lor x_i^{-1}\}|$, that is, the number of occurrences of the generator $x_i$ (and its inverse) in the word $w$. The feature vector $f_0$ is then
$$
f_0=\langle C(w,x_1),\ldots,C(w,x_N)\rangle.
$$
\item $f_1$ (\emph{Weighted Generator Count}) - The feature vector $f_1$ is the same as $f_0$ above, except it is weighted by the word length $|w|$:
$$
f_1=\frac{1}{|w|}\langle C(w,x_1),\ldots,C(w,x_N)\rangle.
$$
\item$f_2$ through $f_7$ (\emph{Counting Subgraphs}) - Let $U_l=\{u_j\in F(X)\mid |u_j|=l\}$, and consider the previously defined counting functions $C(w,xu_{lj}y)$, with $u_{lj}\in U_l$ and $x,y\in Y$ such that $xu_{lj}y$ is a geodesic word representing the element $\phi(xu_{lj}y)$. The features below represent counting subgraphs as described above, and for each length there is a weighted and non-weighted variant:
$$
\begin{array}{lr}
f_2=&\langle C(w,xu_{1j}y)\mid x,y\in Y; u_j\in U_1\rangle\\
f_3=&\frac{1}{|w|}\langle C(w,xu_{1j}y)\mid x,y\in Y; u_j\in U_1\rangle\\
f_4=&\langle C(w,xu_{2j}y)\mid x,y\in Y; u_j\in U_2\rangle\\
f_5=&\frac{1}{|w|}\langle C(w,xu_{2j}y)\mid x,y\in Y;u_j\in U_2\rangle\\
f_6=&\langle C(w,xu_{3j}y)\mid x,y\in Y; u_j\in U_3\rangle\\
f_7=&\frac{1}{|w|}\langle C(w,xu_{3j}y)\mid x,y\in Y; u_j\in U_3\rangle\\
\end{array}
$$
\end{itemize}
\section{Model Selection and Other System Components}
\subsection{Model Selection}
\label{models}
\paragraph{\nopunct}
In the context of machine learning, a \emph{model} or \emph{learning algorithm} is the means by which a set of training inputs can be used to predict the output on future, unseen inputs. The choice of a learning algorithm is informed by the type and structure of the training data, such as whether the data is discrete or continuous, or is comprised of feature vectors like those described above. A particular learning algorithm in turn determines the \emph{hypothesis space}; the set of functions that can be learned from the data. The class of available learning algorithms that have been developed is too numerous to describe here. Instead, we will focus on a set of models that will be applied to group-theoretic problems: decision trees, random forests, and $N$-tuple networks.
\subsubsection{Decision Trees and Random Forests}
\paragraph{\nopunct}
\emph{Decision tree learning} is a model that utilizes a tree structure to encode the learned function. Trees can be used for either classification or regression analysis; we will focus on those used for classification, and in particular binary classification trees, where each node can have a maximum of two children. Every node in the decision tree corresponds to a unique partition of the measurement space. Leaf nodes correspond to the assignment of a particular class, while at internal nodes a test is performed to distinguish between data points. This distinction is encoded in the node's children and further partitions the measurement space. Trees can distinguish by feature, by combining features via a discriminant function (such as a linear discriminant), or by other means.\\
\paragraph{\nopunct}
There are a number of tests available that can be used to partition the space at each internal node. \emph{Gini impurity} measures the frequency at which the remaining data points would be misclassified, with the best split being that which minimizes this impurity. With \emph{information gain}, the entropy of the parent node and the remaining data are calculated, and a partition is chosen that reduces the entropy the most, i.e., that which maximizes the information that can be obtained from the partition. The probability of misclassification can also be used.\\
\paragraph{\nopunct}
As with many learning algorithms, decision trees are prone to overfitting. As decision trees do not have a fixed size representation (i.e., they are considered a non-parametric model), the tree-making algorithm can create large trees or ones with complex branching that do not generalize well. This can be combatted by \emph{pruning}, whereby a subtree of the learned decision tree is replaced with a leaf node whose class is the most common one of the data points contained in the pruned subtree.  Pruning can be performed by testing a subtree's classification performance against a separate data set, and keeping the subtree if improves the performance of the classifier. Pruning can also be achieved by employing a statistical significance test such as the $\chi^2$ test that can determine if a subtree results in a meaningful split in the data, and pruning if the subtree does not meet some threshold of significance. Yet another technique is to limit the depth to which the tree can grow during the training process, a form of \emph{pre-pruning}.\\
\paragraph{\nopunct}
Random forest classifiers \cite{breiman2001random} are an example of an \emph{ensemble} method, in which multiple classifiers are combined into a single classifier. As the name implies, random forests are comprised of several decision trees that are constructed from a random sampling of the training set. Additionally, the best split at each node in a particular tree is determined not by the single best component of the feature vector, but instead by choosing the best feature among a randomly sampled subset of the feature vector's components. Using multiple trees trained on the sampled training set reduces overfitting, while using subsets of the feature vector for choosing partitions in the measurement space reduces variance. Once trained, the classification of a new sample can be determined by either averaging the classifications of each tree in the forest, or by having each tree vote for the sample's class and assigning the sample to the class with a plurality of the votes.
 
\subsubsection{$N$-tuple Neural Networks}
\paragraph{\nopunct}
Another model that we will investigate is the \emph{$N$-tuple neural network}. $N$-tuple classifiers were introduced in 1959 by Bledsoe and Browning \cite{bledsoe1959pattern} as a means of performing printed character recognition using specialized table lookup hardware. In its original implementation, the positions of a binary number $s$ are sampled using a total of $M$ random patterns of size $N$ for each class $C$. The number of classes is dependent upon the application at hand. For instance, if the $N$-tuple classifier was used to classify single digits (0-9), then the number of classes $C$ would be $10$. The sampled positions produce another binary number $b_m$. These samplings are stored in tables $T_{mc}$ (one for each class $c$ and random sample $m$), and the value of the table entry $T_{mc}(b_m)$ is the total number of times a sample $s$ of class $c$ was mapped by $m$ to the binary number $b_m$. A new sample $s'$ is classified by choosing the class for which the sum of $T_{mc}(b'_m)$ is maximal. If no maximum exists, the classifier does not choose a class and is said to reserve decision.\\
\paragraph{\nopunct}
With the advent of radial basis function networks and other forms of artificial neural networks in the late 1980s and early 1990s, $N$-tuple classifiers were revisited, being recast as a type of weightless, single-layer artificial neural network (e.g., the ``single layer lookup perceptrons" of Tattersall et al.\ \cite{tattersall1991single}). In a series of papers, Allinson and Ko\l{}cz extended $N$-tuple neural networks further, developing a binary encoding scheme based on CMAC (cerebellar model arithmetic computer, an older form of artificial neural network) and Gray codes \cite{kolcz1994application}, as well as using NTNNs for regression analysis \cite{kolcz1996n} instead of classification.
In \cite{rohwer1998theoretical}, Rohwer performed a series of experiments on standard pattern recognition databases, finding that for most data sets the best results were achieved with $N$-tuples of size 8 (i.e., a pattern of length $N=8$) and the total number of patterns $M$ around 1000.\\
\paragraph{\nopunct}
NTNNs can be generalized beyond classifying binary data by utilizing the framework of relational algebra. In this rendition, each data point $s$ is transformed into a feature vector $x$ of fixed length $N$.  Let $J_1,\ldots, J_M$ be index sets or \emph{patterns} of uniform length $P$, that is, subsets of the set $\{0,\ldots,N-1\}$ that represent the indices at which to sample $x$. For each index set $J_m$ and class $c\in C$ we form the table $T_{mc}$. Let $\pi_{mc}(x)$ be the projection operation that samples the feature vector $x$ of class $c$ at the indices given by $J_m$. The resultant value $k= \pi_{mc}(x)$ is computed and the table $T_{mc}$ is updated: $T_{mc}(k)\leftarrow T_{mc}(k)+1$.\\
\paragraph{\nopunct}
Let us illustrate the above operations with a concrete example. Consider an NTNN with the parameters $N=5$, $M=2$, and $P=3$. The table below represents the NTNN table entries for class $0$ after training on the first 3 samples:
\begin{table}[H]
\centering
\begin{tabular}{|c|c|c|rl|}
\hline
$s$&$J_0$&$\pi_{00}$&\multicolumn{2}{|c|}{$T_{00}$}\\
\hline
$\langle -4,-1,5,2,3\rangle$&\multirow{3}{*}{$(0,2,4)$}&$(-4,5,3)$&$(-4,5,3)\mapsto$&2\\
$\langle -4,-7,5,2,3\rangle$&&$(-4,5,3)$&$(-2,6,1)\mapsto$&1\\
$\langle -2,-1,6,3,1\rangle$&&$(-2,6,1)$&&\\
\hline
$s$&$J_1$&$\pi_{10}$&\multicolumn{2}{|c|}{$T_{10}$}\\
\hline
$\langle -4,-1,5,2,3\rangle$&\multirow{3}{*}{$(1,2,3)$}&$(-1,5,2)$&$(-1,5,2)\mapsto$&1\\
$\langle -4,-7,5,2,3\rangle$&&$(-7,5,2)$&$(-7,5,2)\mapsto$&1\\
$\langle -2,-1,6,3,1\rangle$&&$(-1,6,3)$&$(-1,6,3)\mapsto$&1\\
\hline
\end{tabular}
\caption{Table Entries of a NTNN Classifier After Training on 3 Samples}
\label{ntnnex}
\end{table}
\paragraph{\nopunct}
The first pattern $J_0$ will sample the first, third, and fifth component of each training sample $s$. The projection operation using pattern $J_0$, denoted $\pi_{00}$, results in the projection $(-4,5,3)$ twice and $(-2,6,1)$ once. The table $T_{00}$ records both projections, as well as the frequency at which these projections have been observed so far. Table $T_{10}$ records the projections and their frequency resulting from the projection operation $\pi_{10}$, which uses pattern $J_1$ instead of $J_0$. Note that $\pi_{10}$ results in three different projections, thus the frequency of each projection in table $T_{10}$ is set to one.\\
\paragraph{\nopunct}
Various classification criteria \cite{haralickntuple} have been devised for use with NTNNs; we include two such criteria below. Note that if there is no unique class that satisfies the criterion, the NTNN reserves decision:\begin{itemize}
\item \emph{Voting Majority} -- Assign $s$ to class $c'$ if  $\sum_{m\in M}T_{mc'}(s)>\sum_{m\in M}T_{mc}(s)$ for all classes $c\neq c'$.
\item \emph{Logarithm Voting Majority} -- Assign $s$ to class $c'$ if  $\sum_{m\in M}\log T_{mc'}(s)>\sum_{m\in M}\log T_{mc}(s)$ for all classes $c\neq c'$.
\end{itemize}

\

Note that after all the counts are created in all the tables, each count may be replaced by the log of the count in accordance with what classification criterion is used in the experiment.\\
\subsection{Data Generation}
\paragraph{\nopunct}
The primary goal of a supervised machine learning system is to predict classes or values on new, unseen members of the data domain. As a consequence, the system's construction requires at least two sets of data - one for training and one for evaluation. When data is scarce the method of \emph{cross-validation} can be employed, whereby the data set is partitioned into separate sets that are used exclusively for training or evaluation. Depending on the learning model chosen, a third data set may be required for use during the optimization process that occurs between the training and evaluation phases. Our approach is to produce three, independently generated sets that are used in each phase of the system's construction. The sets used for training, optimization, and verification are respectively referred to as $S_i$, $S_o$, and $S_v$.\\
\subsection{Evaluation}
\label{mlevaluation}
\paragraph{\nopunct}
There are various means by which the performance of a machine learning system can be evaluated. Regardless of which method is chosen, it is imperative that an independent data set (e.g., the verification set $S_v$ from the previous section) is used to measure the system's performance.\\
\paragraph{\nopunct}
In classification, we are primarily concerned with the \emph{accuracy} $\mathcal{A}$ of a model $\mathcal{M}(f)$, trained with respect to the feature vector $f$, over the verification set $S_v$:
$$
\mathcal{A}(\mathcal{M}(f),S_v)=\frac{|\textrm{True Positives}(S_v)| + |\textrm{True Negatives}(S_v)|}{|S_v|}.
$$
\section{A Machine Learning Approach to the Conjugacy Decision Problem}
\paragraph{\nopunct}
Recall that the \emph{conjugacy decision problem} for a group $G$ is to determine for any $u,v\in G$ if $u$ is conjugate to $v$. With respect to computability, the conjugacy decision problem is in fact \emph{two} problems - each concerned with determining positive or negative solutions exclusively. The positive conjugacy decision problem in any recursively presented group is computable, as for any element in the group its conjugates can be recursively enumerated \cite{Myasnikov-book}. The negative solution is not guaranteed to be computable for non-finite groups. There are classes of groups for which both parts of the conjugacy decision problem are computable, including finitely generated polycyclic and metabelian groups.\\
\paragraph{\nopunct}
Computability, however, does not imply efficiency. The efficient algorithms that do exist are often restricted in some sense, such as answering only one part of the decision problem or being solely applicable to a specific class of groups. For instance, a polynomial algorithm exists \cite{lysenok2010conjugacy} for the full conjugacy decision problem in the Grigorchuk groups, and a linear algorithm exists \cite{epstein2006linearity} for word-hyperbolic groups. In non-cyclic finitely generated groups of infinite abelianization, a linear algorithm was found \cite{kapovich2003generic} that can be used to solve the negative conjugacy decision problem \emph{generically}, i.e., for ``most inputs" (a rigorous formulation of which is introduced in the same paper  \cite{kapovich2003generic}).\\
\paragraph{\nopunct}
Given the limitations of existing algorithms, we turn to the framework of the previous section for a machine learning solution. In the sections below we outline how we adapt the general framework to constructing machine learning systems for the conjugacy decision problem. We use the aforementioned supervised learning methods to train classifiers for several classes of finitely presented group. These classifiers can determine whether a pair of elements from their respective groups are conjugate or not, and do so with very high accuracy.
\subsection{Feature Extraction}
\paragraph{\nopunct}
Given a group $G$ with an efficiently calculable normal form and words $u,v\in G$, we concatenate (denoted by $\parallel$) the unit feature vectors $n_0$ and $n_1$ from section \ref{mlfv} to create two derived feature vectors for the conjugacy decision problem:
$$
\begin{array}{l}
c_0=\langle n_0(u)\parallel n_0(v)\rangle\\
c_1=\langle n_1(u)\parallel n_1(v)\rangle\\
\end{array}
$$
The default feature vector used with tree-based classifiers is $c_1$ (weighted normal forms), while for NTNNs it is $c_0$ (unweighted normal forms). 
Additional feature vectors and normal forms that are used for a particular group are included in that group's experimental results section.
\subsection{Model Selection}
\paragraph{\nopunct}
For the conjugacy decision problem we can utilize all three classifiers defined previously: decision trees, random forests, and NTNNs.\\
\paragraph{\nopunct}
In constructing our decision tree and random forests classifiers we used the implementations included in Scikit-learn \cite{scikit-learn} (specifically the  \texttt{DecisionTreeClassifier}  and \texttt{RandomForestClassifier}).  We trained an instance of each classifier for each group on its respective training set $S_i$. Each leaf node was required to contain at least one sample, and a split could only occur if there were at least two samples at that node. The accuracy of each classifier was calculated by classifying the data in each group's verification set $S_v$. Both Gini impurity and information gain were used to determine the best split. For the depth limit, we tested not having a limit, as well as limiting the number of levels to $\log_2 S_i-1$.  For the random forest classifier, the number of trees in the forest was 10, and the size of the random subset of features was $\sqrt{|c_1|}$, i.e., the square root of the length of the feature vector.\\
\paragraph{\nopunct}
Our version of the NTNN classifier was implemented in Python. We trained a number of NTNN classifiers for each group with different combinations of $M$ patterns of size $P$. The total number of patterns $M$ was varied, taking on values from the set $\{10,20,30,50,100\}$.  The initial size of the patterns was set to 3, and where applicable, sizes in the range $[3,5]$ were tested. We used both the ``Voting Majority" and ``Logarithm Voting Majority" criteria for classification in our tests.\\
\paragraph{\nopunct}
For a given feature vector of dimension $N$, the total number of patterns of size $P\leq N$ is $\binom{N}{P}$. When initializing a NTNN classifier that uses pattern sets of size $M$ (a set of $M$ patterns of size $P$), the list of $\binom{N}{P}$ patterns is generated, and a separate permutation of each list is kept for each class $C$. In our case, $C=2$, as we are performing binary classification, i.e., determining if the given feature vector represent a conjugate or non-conjugate pair of elements. Before training the NTNN on the set $S_i$, each class is assigned the first $M$ patterns from its pattern list.\\
\paragraph{\nopunct}
As the accuracy of the initial random selection of patterns varies considerably, a random restart was implemented, in which a new NTNN was initiated with a new random permutation of all possible patterns. Each NTNN's performance was tested against the set $S_o$, with the NTNN proceeding to the optimization stage only when its accuracy was greater than the \emph{starting threshold} $\theta_{\alpha}$, which was set to 60\%.\\
\paragraph{\nopunct}
During the optimization phase, the algorithm alternates between each classes' list of patterns in choosing the next test pattern. Each pattern in $M$ is swapped out with the test pattern, and the NTNN is evaluated against the optimization set $S_o$. The algorithm keeps track of the pattern $m$ whose replacement with the test pattern improves accuracy the most over all $m\in M$, and makes that pattern swap permanent if a new best accuracy is achieved. The algorithm will continue this process until all patterns have been exhausted or the \emph{goal threshold} $\theta_{\omega}$ is reached, which was set to 97\%. The NTNN classifier and the current location in the pattern list are then saved, so that optimization can be continued at a later time if desired.
\subsection{Data Generation}
\label{datagen}
\paragraph{\nopunct}
In solving the CDP, one instance of each of the three data sets was generated for each group. Each data set consists of 20,000 geodesic word pairs, 10,000 of which represent a pair of conjugate elements, while the other  10,000 represent pairs of non-conjugate elements. These data sets are generated via the following procedures:
\begin{enumerate}
\item \label{ncgen}\emph{Random Non-Conjugate Word Pairs in Normal Form} - For each $n\in[5,1004]$ we generate two words  $u,v$ representing elements $x,y\in G$ respectively, with $|u|=|v|=n$. A word $w$ is generated uniformly and randomly by starting with the identity element $w=1_G$, then selecting a generator (or inverse) $g$ from $X\cup X^{1}$ and performing the product $w=w\cdot g$. The element is then converted into its normal form $w'$, and the length $|w'|$ is computed. Additional products are computed until $|w'|=n$.\\
\paragraph{\nopunct}
After generating each $u,v$ pair, an additional step is required to verify that $u$ is not conjugate $v$. Using the method from \cite{kapovich2003generic}, we construct the derived (or commutator) subgroup of $G$, denoted $[G,G]$, and an epimorphism $\psi:G\rightarrow G/[G,G]$. If we let $\phi$ be the canonical epimorphism, we then look at the images $\psi(\phi(u))$ and $\psi(\phi(v))$ and reject the pair if they map to the same representative in the quotient $G/[G,G]$ or if either maps to $1_{G/[G,G]}$. For $\psi(\phi(u))=\psi(\phi(v))$ if and only if $x\bar{G}=y\bar{G}$ for some coset $\bar{G}$, and this is the case when $y=gxg^{-1}$  for $\phi(u)=x$, $\phi(v)=y$, and some $g\in G$.\\
\paragraph{\nopunct}
This two-step process is repeated until 10 non-conjugate pairs are generated for each $n$.
\item \label{cgen}\emph{Random Conjugate Word Pairs in Normal Form} - For $n\in[5,1004]$ we generate a pair of words $u,t$  representing elements $x,z\in G$ respectively, with $|u|=|t|=n$. Each word $u,t$ is generated uniformly and randomly as above. After $u$ and $t$ are generated, the word $v=u^t$ is formed, and the tuple $(u,v)$ is added to the dataset (z is discarded). This process is repeated 10 times for each $n$.\\
\paragraph{}
The above procedures were implemented in the computer algebra system \texttt{GAP} \cite{gap}.
\end{enumerate}
\paragraph{}
To summarize, the data generation process above for each group $G$ produces a collection of data sets $D_0(G)=\{S_i,S_o,S_v\}$, each with 20,000 pairs of words, with the following properties for each class (conjugate pairs vs. non-conjugate pairs):
\begin{itemize}
\item 10,000 pairs of conjugate words $(u,v)$, with $|u|=l$ for $l\in[5,1004]$ and $|v|$ varying;
\item 10,000 pairs of non-conjugate words $(u,v)$, with $|u|=|v|=l$.
\end{itemize}
These data collections $D_0(G)$ were used to evaluate the parameters of each machine learning model for each group $G$. However, the very positive experimental results reported in section \ref{expresults} may be due to the particular characteristics of the lengths of the words within each class. Therefore, we generated additional dataset collections of increasing generality in order to allay these concerns.\\
\paragraph{}
We first generated large sets of words in normal form for each group by randomly and uniformly selecting lengths in the range $[5,1004]$ and using the word generation procedure described in (\ref{ncgen}) above (without testing for non-conjugacy). This process was continued until we generate a dataset $P(G)$ containing a minimum of 250,000 unique words in normal form for each group $G$, which was a sufficient number of unique words to produce the additional data collections described below.\\
\paragraph{}
The data collection $D_1(G)$ keeps the same subset of conjugate word pairs as $D_0(G)$ for each data set $S_k$, but replaced the non-conjugate word pairs by drawing two words $u,v$ uniformly without replacement from $P(G)$, and adding the pair $(u,v)$ to the dataset if their non-conjugacy was confirmed using procedure described in (\ref{cgen}) above. Each data set $S_k\in D_1(G)$ contains 10,000 conjugate word pairs $(u,v=u^t)$ with $|u|=|t|=l$ for $l\in[5,1004]$, and $|v|$ varied; and 10,000 non-conjugate pairs $(u,v)$ with $|u|=m$,$|v|=n$ for $m,n\in[5,1004]$.\\
\paragraph{}
The data collection $D_2(G)$ keeps the same subset of non-conjugate word pairs as $D_1(G)$ for each data set $S_k$, but replaced the conjugate pairs by drawing two words $u,t$ uniformly without replacement from $P(G)$, creating the word $v=u^t$, and adding the pair $(u,v)$ to the dataset. Each data set $S_k\in D_2(G)$ contains 10,000 conjugate word pairs $(u,v=u^t)$ with $|u|=l$,$|t|=p$ for $l,p\in[5,1004]$, and $|v|$ varied; and 10,000 non-conjugate pairs $(u,v)$ with $|u|=m$,$|v|=n$ for $m,n\in[5,1004]$.\\
\paragraph{}
In the data collections $D_2(G)$ we have ensured that each word regardless of class can vary in length. However, it is still the case that non-conjugate words have lengths strictly within the range $[5,1004]$, while the length of conjugate words have a much larger range. To remove this restriction we created a third data collection $D_3(G)$ for each group. $D_3(G)$ contains the same subset of conjugate pairs as $D_2(G)$. For the non-conjugate pairs, we first determined for each dataset $S_K\in D_2(G)$ the minimum and maximum lengths of conjugate words, and produced non-conjugate word pairs $(u,v)$ with $|u|\in[5,1004]$ as before, but the length of $v$ was allowed to vary within the minimum and maximum lengths of conjugate words. This ensures that length alone can not be used to determine if a pair is conjugate or not.\\
\paragraph{}
For convenience and clarity we have summarized the differences in these data collections in table below. Let $\min_{ik}$ ($\max_{ik}$) correspond to the minimum (maximum) word length in a conjugate pair for the data set $S_k\in D_i(G)$:
\begin{table}[H]
\centering
\begin{tabular}{|c| rr | rr |}
\hline
\textbf{Collection}&\multicolumn{2}{c|}{\textbf{Conjugate Pair} $(u,v=u^t)$}&\multicolumn{2}{c|}{\textbf{Non-Conjugate Pair} $(u,v)$}\\
\hline
$D_0$ & $|u|=|t|=l$;&$l\in[5,1004]$ & $|u|=|v|=m$;& $m\in[5,1004]$\\
\hline
$D_1$ & $|u|=|t|=l$;&$l\in[5,1004]$ & $|u|=m,|v|=n$;& $m,n\in[5,1004]$\\
\hline
$D_2$ & $|u|=l$,$|t|=p$;&$l,p\in[5,1004]$ & $|u|=m,|v|=n$; &$m,n\in[5,1004]$\\
\hline
\multirow{2}{*}{$D_3$}	& \multirow{2}{*}{$|u|=l$,$|t|=p$;} &\multirow{2}{*}{$l,p\in[5,1004]$}	& $|u|=m,|v|=n$; &$m\in[5,1004]$\\
					&&							& \multicolumn{2}{c|}{$n\in[\min_{2k},\max_{2k}]$}\\
\hline
\end{tabular}
\caption{Word Length Differences in Each Data Collection and Class}
\label{table:datacolls}
\end{table}

\subsection{Test Groups}
\paragraph{\nopunct}
Having specified how we will apply our machine learning system to the conjugacy decision problem, we now turn our attention to the groups in which we will evaluate the system's performance. We tested six groups with various algebraic properties:
\begin{itemize}
\item three non-virtually nilpotent groups: the groups $\mathcal{O}\rtimes U_H$;
\item two non-polycyclic metabelian groups: BS(1,2) and GMBS(2,3);
\item and SL$(2,\Z)$, a non-solvable linear group.
\end{itemize}
The specific representations and feature vectors used for each group are outlined below.
\subsubsection{The Baumslag-Solitar Group BS(1,2)}
The Baumslag-Solitar groups are a well-known class of one-relator groups. We will consider the Baumslag-Solitar group BS$(1,2)$, given by the following presentation:
$$
BS(1,2)=\langle a, b, \mid bab^{-1}a^{-2}\rangle.
$$
\paragraph{\nopunct}
Note that the conjugacy decision problem over BS(1,2) resides in the complexity class \textsc{TC}$^0$ \cite{weiss2016logspace}, where \textsc{TC}$^0$ is the class of constant-depth arithmetic circuits using \textsc{AND}, \textsc{OR}, \textsc{NOT}, and majority gates.\\
\paragraph{\nopunct}
Elements in BS(1,2) can be uniquely written in the following normal form:
$$
n_0=b^{-e_1}a^{e_2}b^{e_3},
$$
with $e_1, e_3\geq 0$, and if $e_1, e_3>0$, then $e_2$ is not divisible by 2. Collection from the left can transform any element of BS(1,2) into this normal form.  The dimension of feature vectors $c_0$ and $c_1$ is 6.\\
\paragraph{\nopunct}
In our initial tests of the NTNN classifier for BS(1,2), the classifier did not perform well using the feature vector $c_0$. It may be that the relatively low dimension of the feature vector ($N=6$) provides insufficient information to the classifier. Therefore, we also tested features vectors $c_2$ and $c_4$, defined as the concatenation of unit vectors $f_2$ and $f_4$ respectively. For BS(1,2), the feature vector $c_2$ is of dimension 48, while $c_4$ has dimension 96. 
\subsubsection{Non-Virtually Nilpotent Polycyclic Groups}
\paragraph{\nopunct}
Polycyclic groups are natural generalizations of cyclic groups. A group $G$ is said to be {\em polycyclic} if it has a subnormal series
$$
G=G_1 \triangleright \cdots \triangleright G_{n+1}=\{1\}
$$
 such that the quotient groups $G_i/G_{i+1}$ are cyclic. This series is called a \emph{polycyclic series}. The {\em Hirsch length} of a polycyclic group $G$ is the number of infinite groups in its polycyclic series.\\
\paragraph{\nopunct}
A sequence of elements $X=[x_1,\ldots,x_n]$ such that $\langle x_iG{i+1}\rangle=G_i/G_{i+1}$ is called a \emph{polycyclic sequence} for $G$. In a polycyclic group $G$ with polycyclic sequence $X$, any element $g$ can be represented uniquely in normal form as a product of powers of the generators of $G$:
$$
g=x_1^{e_1}\cdots x_n^{e_n},
$$
with $e_i\in\Z$. The sequence $(e_1,\ldots,e_n)$ is called the \emph{exponent vector} of $g$ with respect to $X$.\\
\paragraph{\nopunct}
Polycyclic groups that are non-virtually nilpotent have exponential growth (\cite{Wolf},\cite{Milnor}) and remain promising candidates for use as platform groups in non-commutative cryptography \cite{gryak2016status}. One method of constructing such groups is through the use of algebraic number fields, as outlined in \cite[\S 8.2.2]{HEO05}.\\
\paragraph{\nopunct}
Given an algebraic number field $F$ with degree $[F:\Q]>1$, one can define two substructures, the \emph{maximal order} $\mathcal{O}(F)$ and the \emph{unit group} $U(F)$. The maximal order is the largest ring of integers of $F$, and consists of those elements in $F$ that are a root of some monic polynomial over $F$ with integral coefficients. The multiplicative group $U(F)$ consists wholly of the non-zero elements of $\mathcal{O}(F)$ that have a multiplicative inverse, i.e., are units. Given these two structures and the aforementioned degree criterion, the semidirect product $\mathcal{O}(F)\rtimes U(F)$ results in an infinite, non-virtually nilpotent polycyclic group.\\
\paragraph{\nopunct}
Below are three specific instances of the $\mathcal{O}(F)\rtimes U(F)$ family of polycyclic groups. The conjugacy search problem over the first two groups was studied in \cite{garber2013analyzing}, in the context of the length-based attack against the AAG key exchange protocol. The groups can be constructed by using the \texttt{MaximalOrderByUnitsPcpGroup} function of the GAP Polycyclic package \cite{Polycyclic2.11}. The function takes a polynomial that is irreducible over $\Q$ (thereby defining a field extension of $\Q$) and returns a group of the form $\mathcal{O}(F)\rtimes U(F)$: 
\begin{itemize}
\item ${\mathcal{O}\rtimes U}_{14}$ - Given the polynomial $f=x^9-7x^3-1$,\\ \texttt{MaximalOrderByUnitsPcpGroup} returns a group of the form $\mathcal{O}(F)\rtimes U(F)$ with a Hirsch length of 14.
\item ${\mathcal{O}\rtimes U}_{16}$ - Given the polynomial $f=x^{11}-x^3-1$,\\ \texttt{MaximalOrderByUnitsPcpGroup} returns a group of the form $\mathcal{O}(F)\rtimes U(F)$ with a Hirsch length of 16.
\item ${\mathcal{O}\rtimes U}_{34}$ - Given the polynomial $f=x^{23}-x^3-1$,\\ \texttt{MaximalOrderByUnitsPcpGroup} returns a group of the form $\mathcal{O}(F)\rtimes U(F)$ with a Hirsch length of 34.
\end{itemize}
\paragraph{\nopunct}
Recall that every polycyclic group has a normal form in terms of the generators in its polycyclic sequence. The feature vector $n_0$ for a polycyclic group element $g$ simply corresponds to the exponent vector of $g$ in normal form.  Thus the feature vectors   $c_0$ and $c_1$ are readily computable for polycyclic group elements. The dimension of these feature vectors for groups of the form $\mathcal{O}(F)\rtimes U(F)$ is $2(H+1)$, where $H$ is the Hirsch length of the group.
\subsubsection{Generalized Metabelian Baumslag-Solitar Groups}
\paragraph{\nopunct}
In \cite{gryak2016conjugacy}, Gryak et al.\ introduced a family of polycyclic and metabelian groups for which the time complexity of the conjugacy search problem is exponentially bounded. Generalized metabelian Baumslag-Solitar groups are a subfamily of these groups, with the group GMBS(2,3) being one such member whose presentation is given below:
$$
\textrm{GMBS}(2,3)=\langle q_1,q_2,b\mid b^{q_1}=b^2, b^{q_2}=b^3, [q_1,q_2]=1\rangle.
$$

\paragraph{\nopunct}
Elements in GMBS(2,3) can be uniquely written in the following normal form:
$$
n_0=q_1^{-e_1}q_2^{-e_2}b^{e_3}q_1^{e_4}q_2^{e_5},
$$
with $e_1,e_2,e_4,e_5\geq 0$, $2\nmid e_3$ if $e_1,e_4>0$, and $3\nmid e_3$ if $e_2,e_5>0$. Collection from the left can transform any element of GMBS(2,3) into this normal form.  The dimension of feature vectors $c_0$ and $c_1$ is 10.
\subsubsection{SL$(2,\Z)$}
Recall that SL$(2,\Z)$ is the set of $2\times2$ integral matrices with determinant 1. This set forms a group under matrix multiplication, and is a discrete subgroup of SL$(2,\R)$. The group was implemented in \texttt{GAP} with a dual representation: for each element $x\in\;$SL$(2,\Z)$  we have a pair $(m,w)$ of the form

$$
m=
\left[
\begin{array}{cc}
a&b\\
c&d\\
\end{array}\right]
,
w=w_1\cdots w_n, w_i\in\{S^{\pm 1},R^{\pm 1}\},
$$

with $a,b,c,d\in \Z$ such that $ad-bc=1$, and $S$ and $R$ corresponding to the matrices below that generate SL$(2,\Z)$:

$$
S=
\left[
\begin{array}{rr}
0&-1\\
1&0\\
\end{array}\right]
,
R=
\left[
\begin{array}{rr}
0&-1\\
1&1\\
\end{array}\right].
$$

\paragraph{\nopunct}
In this formulation, SL$(2,\Z)$  is an amalgamated free product given by the presentation
$$
\text{SL}(2,\Z)\cong\langle S,R\mid S^4=1,S^2=R^3\rangle.
$$
\paragraph{\nopunct}
These generators and attendant presentation were chosen so that a confluent rewriting system could be constructed in \texttt{GAP} via the Knuth-Bendix algorithm. This rewriting system enables us to reduce any word to shortlex normal form efficiently.\\
\paragraph{\nopunct}
In generating the data sets, the length of an element $x$ was taken to be the length of the word representation of the element, i.e., $|x|=|w|$, as suggested in \cite{shpilrainmatrix2017}. When the matrix form $m$ of the element $x$ is needed, the norm of the matrix, $\lVert m\lVert$, can be used as a length measure. We utilized the \emph{Frobenius norm}, which is calculated as
$$
\lVert m\lVert=\sqrt{a^2+b^2+c^2+d^2}.
$$
\label{sl2znf}
\paragraph{\nopunct}
Given that there are two representations for each element, there are multiple normal forms that can be considered. Let $x=(m,w)\in\;$SL$(2,\Z)$. The \emph{matrix normal form} is simply the ``flattened" matrix, i.e., a vector in $\Z^4$ :
$$
f_m=\langle a,b,c,d\rangle.
$$ 
\paragraph{\nopunct}
For the tree-based classifiers we used the normalized matrix normal form as the feature vector, i.e., for a word $u=(m_u,w_u)$ we have 
$$
f_{m}=\frac{1}{\lVert m_u\lVert}\langle a,b,c,d\rangle,
$$
and for a pair of words $u,v$ with respective matrix representations $m_u,m_v$, we concatenate the two unit feature vectors together to form a feature vector for the conjugacy decision problem:
$$
c_m=\langle f_m(m_u)\parallel f_m(m_v)\rangle
$$
\paragraph{\nopunct}
For the NTNN classifier we are looking to train on discretely valued data, thus the previous feature vector of normalized matrix entries is not applicable. Attempting to use the \emph{unnormalized} matrix normal form would be a poor choice, as the frequency distribution of the integral values that comprise the matrix entries is highly skewed. Thus in lieu of the matrix representation for an element we will use its word representation. However, the shortlex normal form for SL$(2,\Z)$ does not have a fixed length. Consequently, we will use the counting subgraph features to transform each reduced word in our formulation of SL$(2,\Z)$ (as an amalgamated free product) to a fixed length feature vector. In particular, we will utilize the feature vector $c_2$ that was used for BS(1,2). Note however that for SL$(2,\Z)$ the dimension of $c_2$ is 40.
\subsection{Experimental Results}
\label{expresults}
\paragraph{\nopunct}
In this section we present the results of the performance of our three classification models on the above test groups. While we tested the all three models over all combinations of their respective parameters, the accuracy results in the tables below represent the single best performing classifier for each group and model. Unless otherwise noted, the feature vector $c_1$ (weighted normal forms) was used for the tree-based classifiers, while $c_0$ (unweighted normal forms) was used for the NTNN classifier.\\
\paragraph{\nopunct}
The overall accuracy of a classifier for a particular group masks differences in how that classifier performs over different samples in the verification set. To elucidate these differences, we present two more granular analyses of accuracy: that with respect to word length and another with respect to class.
\subsubsection{Decision Tree and Random Forests}
\paragraph{\nopunct}
Tables \ref{dtperf} and \ref{rfperf} below respectively display the accuracy of the best performing decision tree and random forest classifiers for each group. For all classifiers in this section the feature vector $c_1$ was used, with the exception of the classifiers for SL$(2,\Z)$, which used the matrix-based vector $c_m$.
\begin{table}[H]
\centering
\begin{tabular}{|c|l|c|}
\hline
\textbf{Group}&\multicolumn{1}{|c|}{\textbf{Method, Split Criterion, Depth}}&\textbf{Accuracy}\\
\hline
BS(1,2)&Decision Tree, Entropy, Depth Limit&92.00\%\\
${\mathcal{O}\rtimes U}_{14}$&Decision Tree, Entropy, Depth Limit&98.49\%\\
${\mathcal{O}\rtimes U}_{16}$&Decision Tree, Entropy, No Depth Limit&97.23\%\\
${\mathcal{O}\rtimes U}_{34}$&Decision Tree, Entropy, Depth Limit&98.47\%\\
GMBS(2,3)&Decision Tree, Gini Impurity, Depth Limit&95.43\%\\
SL$(2,\Z)$&Decision Tree, Entropy, No Depth Limit&96.26\%\\
\hline
\end{tabular}
\caption{Best Performing Decision Tree Classifiers for All Groups}
\label{dtperf}
\end{table}
In optimizing the performance of the decision tree-based classifiers, different combinations of tree depth limits and splitting criteria were considered. For nearly all test groups, using information gain (equivalently, greatest reduction in entropy) resulted in the most accurate classifier. Only for GMBS(2,3) did using Gini impurity result in a higher accuracy, and only by .1\%. \\
\paragraph{\nopunct}
\begin{table}[H]
\centering
\begin{tabular}{|c|l|c|}
\hline
\textbf{Group}&\multicolumn{1}{|c|}{\textbf{Method, Split Criterion, Depth}}&\textbf{Accuracy}\\
\hline
BS(1,2)&Random Forest, Entropy, No Depth Limit & 93.64\%\\
${\mathcal{O}\rtimes U}_{14}$&Random Forest, Entropy, No Depth Limit & 98.69\%\\
${\mathcal{O}\rtimes U}_{16}$&Random Forest, Entropy, Depth Limit& 98.19\%\\
${\mathcal{O}\rtimes U}_{34}$&Random Forest, Entropy, No Depth Limit & 98.89\%\\
GMBS(2,3)&Random Forest, Entropy, No Depth Limit & 96.49\%\\
SL$(2,\Z)$&Random Forest, Entropy, No Depth Limit & 97.47\%\\
\hline
\end{tabular}
\caption{Best Performing Random Forest Classifiers for All Groups}
\label{rfperf}
\end{table}
For all groups tested, the random forest classifier performed better than a single decision tree, and again information gain resulted in the most accurate classification. Limiting the depth of the tree (or trees in the case of random forests)  to $\log_2 N -1$, where $N$ is the total number of samples, slightly improved the results when using Gini impurity as the splitting criterion, but did not do so when using information gain.\\
\paragraph{\nopunct}
The generalization error for random forests approaches zero as additional trees are included in the forest. For example, the table below lists the classification accuracy for random forest classifiers on the group ${\mathcal{O}\rtimes U}_{34}$ with different numbers of trees. Note the diminishing marginal increases in accuracy as the number of trees increases (as this is a stochastic process, increases in accuracy are not strictly monotonic):
\begin{table}[H]
\centering
\begin{tabular}{|c|c|}
\hline
\textbf{\# Trees}&\textbf{Accuracy}\\
\hline
10&98.89\%\\
15&99.17\%\\
20&99.07\%\\
30&99.20\%\\
50&99.31\%\\
100&99.39\%\\
200&99.41\%\\
\hline
\end{tabular}
\caption{Accuracy of Random Forest Classifiers for ${\mathcal{O}\rtimes U}_{34}$ with Increasingly Large Forests}
\end{table}

\subsubsection{NTNNs}
\paragraph{\nopunct}
Table \ref{ntnnperf} below displays the accuracy of the best performing NTNN classifiers for each group. The table indicates the total number of patterns $M$ and pattern size $P$ used with each classifier. The decision criterion used for each classifier is recorded in parentheses to the right of the accuracy entry, with $(\Sigma)$ indicating that voting majority was used, while $(\log)$ indicating that logarithm voting majority was used instead. 
\begin{table}[H]
\centering
\begin{tabular}{|c|c|c|l|}
\hline
\textbf{Group}&$\mathbf{M}$&$\mathbf{P}$&\textbf{Accuracy}\\
\hline
BS(1,2) & 30 & 4  & 92.41\% (log)\\
${\mathcal{O}\rtimes U}_{14}$ & 20 & 3  & 98.77\% (log)\\
${\mathcal{O}\rtimes U}_{16}$& 20 & 5 & 98.46\% ($\Sigma$)\\
${\mathcal{O}\rtimes U}_{34}$ & 100 & 3 & 99.50\% (log)\\
GMBS(2,3) & 30 & 4 & 96.13\% ($\Sigma$)\\
SL$(2,\Z)$ & 50 & 4 & 99.81\% (log)\\
\hline
\end{tabular}
\caption{Best Performing NTNN Classifiers for All Groups}
\label{ntnnperf}
\end{table}
\paragraph{\nopunct}
For the group BS(1,2), the use of the feature vector $c_2$ produced a marked improvement in accuracy as compared to $c_0$. The accuracy of the NTNN classifier depicted in Table \ref{ntnnperf} above uses $c_2$. We ran the full array of tests over all $(M,P)$ pairs for $c_2$, but ran only three additional tests using the $c_4$ feature vector, as in these initial tests we did not see any improvement in the performance of the NTNN classifier as compared to using $c_2$.\\
\paragraph{\nopunct}
The NTNN classifier for SL$(2,\Z)$ depicted in Table \ref{ntnnperf} was not only the best performing NTNN classifier but the best performing classifier for any model and group tested. The feature vector used for this classifier was also $c_2$, which is of dimension 40 as compared to 48 for BS(1,2).\\
\paragraph{\nopunct}
The classifiers for the remaining groups utilized the feature vector $c_0$. All NTNN classifiers for the ${\mathcal{O}\rtimes U}_{H}$ groups had accuracies above 98\%. The classifier for GMBS(2,3), with an accuracy above 96\%, performed better than that for BS(1,2) but worse than those for the other test groups.
\subsubsection{Accuracy with Respect to Word Length}
\paragraph{\nopunct}
The length of a word with respect to a generating set is a crucial measurement throughout group theory. Word length, rather than bit length, is the standard input size parameter for group-theoretic algorithms. The growth rate of a group, which depends on word length, can determine algebraic properties such as nilpotency. Recall that in non-commutatively cryptography, the word length corresponds to key size. Thus, it is important to consider how well our system performs with respect to the word length.\\
\paragraph{\nopunct}
In analyzing the performance of our classifiers, we looked for a length threshold $L$ that would provide the greatest difference in accuracy between words below and above this demarcation. To calculate $L$ for each group, we first calculated the accuracy of the best performing NTNN classifier for each length and class over the verification dataset $S_v$. We then determined the inflection points in this data via second order finite differences. The threshold $L$ was then set to the length that resulted in the greatest difference in accuracy. The results for each class and group are listed in the table below:
\begin{table}[H]
\centering
\begin{tabular}{|c|c|c|c|c|}
\cline{4-5}
\multicolumn{3}{c}{}&\multicolumn{2}{|c|}{Accuracy}\\
\hline
Group&Class&$L$&$|w|<L$&$|w|\geq L$\\
\hline
\multirow{2}{*}{BS(1,2)}&Conjugate&16&30.00\%&88.82\%\\
&Non-Conjugate&14&84.44\%&96.75\%\\
\hline
\multirow{2}{*}{ ${\mathcal{O}\rtimes U}_{14}$}&Conjugate&10&94.00\%&99.98\%\\
&Non-Conjugate&11&80.00\%&97.69\%\\
\hline
\multirow{2}{*}{ ${\mathcal{O}\rtimes U}_{16}$}&Conjugate&7&95.00\%&99.51\%\\
&Non-Conjugate&21&86.25\%&97.59\%\\
\hline
\multirow{2}{*}{ ${\mathcal{O}\rtimes U}_{34}$}&Conjugate&7&55.00\%&99.23\%\\
&Non-Conjugate&36&97.74\%&99.93\%\\
\hline
\multirow{2}{*}{GMBS(2,3)}&Conjugate&17&88.33\%&97.48\%\\
&Non-Conjugate&9&100\%&94.86\%\\
\hline
\multirow{2}{*}{SL$(2,\Z)$}&Conjugate&17&90.83\%&99.98\%\\
&Non-Conjugate&7&90.00\%&99.77\%\\
\hline
\end{tabular}
\caption{Accuracy with Respect to Word Length and Class for Tested Groups}
\end{table}
\paragraph{\nopunct}
From the above table one can readily observe that classification is more accurate on longer words than shorter ones, with the only exception being the non-conjugate elements of GMBS(2,3). For BS(1,2) the classifier performed very poorly on short conjugate pairs. For the non-virtually nilpotent polycyclic groups, conjugate pair accuracy was over 99\% for words of length greater than 10, while for non-conjugate pairs the length threshold required to achieve this performance level increased along with Hirsch length. The SL$(2,\Z)$ classifier performed very well in both classes with words greater than 17 in length.
\subsubsection{Accuracy with Respect to Class}
\paragraph{\nopunct}
By examining the confusion matrices for the best classifier for each group, we can observe the accuracy for each class of elements in our data set. The accuracies depicted in the table below are for the best performing NTNN classifier for each group, which, with the exception of BS(1,2), is the best performing classifier for all groups. All classifiers achieved higher accuracy on the class of conjugate elements than the class of non-conjugate elements, with the exception of BS(1,2), which in addition had the lowest accuracy results of all groups tested.
\begin{table}[H]
\centering
\begin{tabular}{|c|c|c|}
\cline{2-3}
\multicolumn{1}{c}{}&\multicolumn{2}{|c|}{Accuracy by Class}\\
\hline
Group&Conjugate&Non-Conjugate\\
\hline
BS(1,2)&88.17\%&96.64\%\\
\hline
${\mathcal{O}\rtimes U}_{14}$&99.95\%&97.58\%\\
\hline
${\mathcal{O}\rtimes U}_{16}$&99.50\%&97.41\%\\
\hline
${\mathcal{O}\rtimes U}_{34}$&99.14\%&99.86\%\\
\hline
GMBS(2,3)&97.37\%&94.88\%\\
\hline
SL$(2,\Z)$&99.87\%&99.75\%\\
\hline
\end{tabular}
\caption{Accuracy by Class for Tested Groups}
\end{table}
\subsubsection{Accuracy on Different Data Collections}
All of the previous experimental results and analysis were performed on models trained and optimized for each group $G$ over the original data collection $D_0(G)$. Recall from section \ref{datagen} that we generated additional data sets $D_1(G),D_2(G)$, and $D_3(G)$ with increasingly varied word lengths to verify that our experimental results were not due to the particular word lengths we chose to use in our original experiments.\\
\paragraph{}
For each data collection new models were trained and optimized for each group. For the new NTNN models we used the number of patterns $M$ and pattern size $P$ of the best performing NTNN model for each group as indicated in Table \ref{ntnnperf}. The table below depicts the results of testing the various models on the different data collections (see Table \ref{table:datacolls} on page \pageref{table:datacolls} for the composition of each data collection). For brevity, we include the accuracy and type of the single best performing model for each group and data collection. Despite the changes in word lengths within each data collection classification accuracy was maintained.
\paragraph{\nopunct}
\begin{table}[H]
\centering
\begin{tabular}{|c| c | c | c | c |}
\cline{2-5}
\multicolumn{1}{}{}&\multicolumn{4}{|c|}{\textbf{Data Collection}}\\
\hline
\textbf{Group}&$\mathbf{D_0}$&$\mathbf{D_1}$&$\mathbf{D_2}$&$\mathbf{D_3}$\\
\hline
BS(1,2) & 93.64\% (F$_e$) & 93.20\% (F$_g$) &95.30\% (F$_e$) & 98.86\% (F$_e$)\\
\hline
${\mathcal{O}\rtimes U}_{14}$ & 98.77\% (N$_l$) & 98.67\% (F$_{ed}$) &  98.38\% (F$_{ed}$) & 99.75\% (N$_l$)\\
\hline
${\mathcal{O}\rtimes U}_{16}$& 98.46\% (N$_s$) & 97.24\% (F$_{ed}$)  & 96.65\% (F$_{ed}$) & 99.11\% (F$_g$)\\
\hline
${\mathcal{O}\rtimes U}_{34}$ & 99.50\% (N$_l$) & 98.72\% (F$_{ed}$) & 98.28\% (F$_{ed}$) & 99.29\% (N$_l$)\\
\hline
GMBS(2,3) & 96.49\% (F$_e$) & 95.22\% (F$_e$) &96.45\% (N$_s$) &99.13\% (F$_{gd}$)\\
\hline
SL$(2,\Z)$ & 99.81\% (N$_l$) & 99.91\% (F$_g$) & 93.89\% (F$_e$) & 97.38\% (F$_g$)\\
\hline
\end{tabular}
\caption{Best Performing Classifiers By Data Collection - The model for each group and data collection is denoted in parentheses next to its accuracy as follows:  Random (F)orest using (e)ntropy or (g)ini coefficient, with (d) indicating a tree depth limit; (N)TNN using voting majority (s) or log voting majority (l).}
\label{genperf}
\end{table}
\section{Conclusion}
\paragraph{\nopunct}
In conclusion, we have shown how the pattern recognition techniques for free groups developed in \cite{haralick2004pattern} can be extended to non-free groups. We demonstrated that the conjugacy decision problem in a variety of groups can be solved with very high accuracy using random forests and $N$-tuple neural networks. For the group BS(1,2) the random forest classifier performed the best, while for all other groups the NTNN classifiers were the most accurate. The NTNN classifier for SL$(2,\Z)$ was the best performing one for any model and group, with an overall accuracy of over 99.8\%.\\
\paragraph{\nopunct}
 As suggested in \cite{haralick2004pattern}, the successful application of pattern recognition techniques to group-theoretic problems can provide experimental evidence for new conjectures in group theory. The decisions made by the decision trees and n-tuple neural network models used in this paper are readily interpretable, thus enabling a computational group theorist to link the results in the model back to their corresponding algebraic inputs.\\
\paragraph{\nopunct}
We have in fact such a potential conjecture at hand. From the high accuracy of the classifiers across the tested groups it is apparent that there is some underlying mathematical relationship with respect to conjugacy that is responsible for the classifiers' performance. We will perform further analysis on the best performing classifiers to tease out what exactly this mathematical relationship is; a forthcoming paper will bring these additional results to light.
\section*{Acknowledgements}
\paragraph{\nopunct}
We would like to thank Benjamin Fine and Vladimir Shpilrain for their helpful suggestions throughout the development of this work. We would also like to thank the reviewers for their constructive questions, comments, and corrections. Delaram Kahrobaei is partially supported by a PSC-CUNY grant from the CUNY Research Foundation, the City Tech Foundation, and ONR (Office of Naval Research) grants N000141210758 and N00014-15-1-2164. 
\bibliographystyle{plain}

\end{document}